\documentstyle[12pt,amssymb]{article}

\begin{document}

\baselineskip15pt

\newtheorem{definition}{Definition $\!\!$}[section]
\newtheorem{prop}[definition]{Proposition $\!\!$}
\newtheorem{lem}[definition]{Lemma $\!\!$}
\newtheorem{corollary}[definition]{Corollary $\!\!$}
\newtheorem{theorem}[definition]{Theorem $\!\!$}
\newtheorem{example}[definition]{\it Example $\!\!$}
\newtheorem{remark}[definition]{Remark $\!\!$}

\newcommand{\nc}[2]{\newcommand{#1}{#2}}
\newcommand{\rnc}[2]{\renewcommand{#1}{#2}}
\nc{\bpr}{\begin{prop}}
\nc{\bth}{\begin{theorem}}
\nc{\ble}{\begin{lem}}
\nc{\bco}{\begin{corollary}}
\nc{\bre}{\begin{remark}}
\nc{\bex}{\begin{example}}
\nc{\bde}{\begin{definition}}
\nc{\ede}{\end{definition}}
\nc{\epr}{\end{prop}}
\nc{\ethe}{\end{theorem}}
\nc{\ele}{\end{lem}}
\nc{\eco}{\end{corollary}}
\nc{\ere}{\hfill\mbox{$\Diamond$}\end{remark}}
\nc{\eex}{\end{example}}
\nc{\epf}{\hfill\mbox{$\Box$}}
\nc{\ot}{\otimes}
\nc{\bsb}{\begin{Sb}}
\nc{\esb}{\end{Sb}}
\nc{\ct}{\mbox{${\cal T}$}}
\nc{\ctb}{\mbox{${\cal T}\sb B$}}
\nc{\bcd}{\[\begin{CD}}
\nc{\ecd}{\end{CD}\]}
\nc{\ba}{\begin{array}}
\nc{\ea}{\end{array}}
\nc{\bea}{\begin{eqnarray}}
\nc{\eea}{\end{eqnarray}}
\nc{\be}{\begin{enumerate}}
\nc{\ee}{\end{enumerate}}
\nc{\beq}{\begin{equation}}
\nc{\eeq}{\end{equation}}
\nc{\bi}{\begin{itemize}}
\nc{\ei}{\end{itemize}}
\nc{\kr}{\mbox{Ker}}
\nc{\te}{\!\ot\!}
\nc{\pf}{\mbox{$P\!\sb F$}}
\nc{\pn}{\mbox{$P\!\sb\nu$}}
\nc{\bmlp}{\mbox{\boldmath$\left(\right.$}}
\nc{\bmrp}{\mbox{\boldmath$\left.\right)$}}
\rnc{\phi}{\mbox{$\varphi$}}
\nc{\LAblp}{\mbox{\LARGE\boldmath$($}}
\nc{\LAbrp}{\mbox{\LARGE\boldmath$)$}}
\nc{\Lblp}{\mbox{\Large\boldmath$($}}
\nc{\Lbrp}{\mbox{\Large\boldmath$)$}}
\nc{\lblp}{\mbox{\large\boldmath$($}}
\nc{\lbrp}{\mbox{\large\boldmath$)$}}
\nc{\blp}{\mbox{\boldmath$($}}
\nc{\brp}{\mbox{\boldmath$)$}}
\nc{\LAlp}{\mbox{\LARGE $($}}
\nc{\LArp}{\mbox{\LARGE $)$}}
\nc{\Llp}{\mbox{\Large $($}}
\nc{\Lrp}{\mbox{\Large $)$}}
\nc{\llp}{\mbox{\large $($}}
\nc{\lrp}{\mbox{\large $)$}}
\nc{\lbc}{\mbox{\Large\boldmath$,$}}
\nc{\lc}{\mbox{\Large$,$}}
\nc{\Lall}{\mbox{\Large$\forall$}}
\nc{\bc}{\mbox{\boldmath$,$}}
\rnc{\epsilon}{\varepsilon}
\rnc{\ker}{\mbox{\em Ker}}
\nc{\ra}{\rightarrow}
\nc{\ci}{\circ}
\nc{\cc}{\!\ci\!}
\nc{\T}{\mbox{\sf T}}
\nc{\can}{\mbox{\em\sf T}\!\sb R}
\nc{\cnl}{$\mbox{\sf T}\!\sb R$}
\nc{\lra}{\longrightarrow}
\nc{\M}{\mbox{Map}}
\rnc{\to}{\mapsto}
\nc{\imp}{\Rightarrow}
\rnc{\iff}{\Leftrightarrow}
\nc{\bmq}{\cite{bmq}}
\nc{\ob}{\mbox{$\Omega\sp{1}\! (\! B)$}}
\nc{\op}{\mbox{$\Omega\sp{1}\! (\! P)$}}
\nc{\oa}{\mbox{$\Omega\sp{1}\! (\! A)$}}
\nc{\inc}{\mbox{$\,\subseteq\;$}}
\nc{\de}{\mbox{$\Delta$}}
\nc{\spp}{\mbox{${\cal S}{\cal P}(P)$}}
\nc{\dr}{\mbox{$\Delta_{R}$}}
\nc{\dsr}{\mbox{$\Delta_{\cal R}$}}
\nc{\m}{\mbox{m}}
\nc{\0}{\sb{(0)}}
\nc{\1}{\sb{(1)}}
\nc{\2}{\sb{(2)}}
\nc{\3}{\sb{(3)}}
\nc{\4}{\sb{(4)}}
\nc{\5}{\sb{(5)}}
\nc{\6}{\sb{(6)}}
\nc{\7}{\sb{(7)}}
\nc{\hsp}{\hspace*}
\nc{\al}{\mbox{$\alpha$}}
\nc{\bet}{\mbox{$\beta$}}
\nc{\ha}{\mbox{$\alpha$}}
\nc{\hb}{\mbox{$\beta$}}
\nc{\hg}{\mbox{$\gamma$}}
\nc{\hd}{\mbox{$\delta$}}
\nc{\he}{\mbox{$\varepsilon$}}
\nc{\hz}{\mbox{$\zeta$}}
\nc{\hs}{\mbox{$\sigma$}}
\nc{\hk}{\mbox{$\kappa$}}
\nc{\hm}{\mbox{$\mu$}}
\nc{\hn}{\mbox{$\nu$}}
\nc{\la}{\mbox{$\lambda$}}
\nc{\hl}{\mbox{$\lambda$}}
\nc{\hG}{\mbox{$\Gamma$}}
\nc{\hD}{\mbox{$\Delta$}}
\nc{\Th}{\mbox{$\Theta$}}
\nc{\ho}{\mbox{$\omega$}}
\nc{\hO}{\mbox{$\Omega$}}
\nc{\hp}{\mbox{$\pi$}}
\nc{\hP}{\mbox{$\Pi$}}
\nc{\bpf}{{\it Proof.~~}}
\nc{\slq}{\mbox{$A(SL\sb q(2))$}}
\nc{\fr}{\mbox{$Fr\llp\! A(SL(2,\IC))\!\lrp$}}
\nc{\slc}{\mbox{$A(SL(2,\IC))$}}
\nc{\af}{\mbox{$A(F)$}}
\rnc{\widetilde}{\tilde}
\nc{\lan}{\langle}
\nc{\ran}{\rangle}
\nc{\phe}{principal homogenous extension}
\nc{\Label}{\label}

\nc{\Section}{\setcounter{definition}{0}\section}
\renewcommand{\theequation}{\thesection.\arabic{equation}}

\newcounter{c}
\renewcommand{\[}{\setcounter{c}{1}$$}
\newcommand{\etyk}[1]{\vspace{-7.4mm}$$\begin{equation}\Label{#1}
\addtocounter{c}{1}}
\renewcommand{\]}{\ifnum \value{c}=1 $$\else \end{equation}\fi}

\newcommand{\dowod}{\noindent{\bf Proof:} }
\newcommand{\Sp}{{\rm Sp}\,}
\newcommand{\Mor}{\mbox{$\rm Mor$}}
\newcommand{\skrA}{{\cal A}}
\newcommand{\Phase}{\mbox{$\rm Phase\,$}}
\newcommand{\id}{{\rm id}}
\newcommand{\diag}{{\rm diag}}
\newcommand{\inv}{{\rm inv}}
\newcommand{\ad}{{\rm ad}}
\newcommand{\poi}{{\rm pt}}
\newcommand{\Dim}{{\rm dim}\,}
\newcommand{\Ker}{{\rm ker}\,}
\newcommand{\Mat}{{\rm Mat}\,}
\newcommand{\Rep}{{\rm Rep}\,}
\newcommand{\Fun}{{\rm Fun}\,}
\newcommand{\Tr}{{\rm Tr}\,}
\newcommand{\supp}{\mbox{$\rm supp$}}
\newcommand{\half}{\frac{1}{2}}

\newcommand{\skrF}{{A}}
\newcommand{\skrD}{{\cal D}}
\newcommand{\skrC}{{\cal C}}

\newcommand{\ttimes}{\mbox{$\hspace{.5mm}\bigcirc\hspace{-4.9mm}
\perp\hspace{1mm}$}}
\newcommand{\Ttimes}{\mbox{$\hspace{.5mm}\bigcirc\hspace{-3.7mm}
\raisebox{-.7mm}{$\top$}\hspace{1mm}$}}
\newcommand{\Cstar}{$^*$-}

\newcommand{\bbr}{{\bf R}}
\newcommand{\bbz}{{\bf Z}}
\newcommand{\Ci}{C_{\infty}}
\newcommand{\Cb}{C_{b}}
\newcommand{\fa}{\forall}
\newcommand{\rrr}{right regular representation}
\newcommand{\wrt}{with respect to}

\newcommand{\qg}{quantum group}
\newcommand{\qgs}{quantum groups}
\newcommand{\cs}{classical space}
\newcommand{\qs}{quantum space}
\newcommand{\po}{Pontryagin}
\newcommand{\ch}{character}
\newcommand{\chs}{characters}

\def\esl{{\mbox{$E\sb{\frak s\frak l (2,{\IC})}$}}}
\def\esu{{\mbox{$E\sb{\frak s\frak u(2)}$}}}
\def\zf{{\mbox{${\IZ}\sb 4$}}}
\def\zt{{\mbox{$2{\IZ}\sb 2$}}}
\def\ox{{\mbox{$\Omega\sp 1\sb{\frak M}X$}}}
\def\oxh{{\mbox{$\Omega\sp 1\sb{\frak M-hor}X$}}}
\def\oxs{{\mbox{$\Omega\sp 1\sb{\frak M-shor}X$}}}
\def\inbar{\,\vrule height1.5ex width.4pt depth0pt}
\def\IC{{\Bbb C}}
\def\IZ{{\Bbb Z}}
\def\IN{{\Bbb N}}
\def\otc{\otimes_{\IC}}
\def\ra{\rightarrow}
\def\ota{\otimes_ A}
\def\otza{\otimes_{ Z(A)}}
\def\otc{\otimes_{\IC}}
\def\h{\rho}
\def\x{\zeta}
\def\th{\theta}
\def\s{\sigma}
\def\t{\tau}
\def\st{\stackrel}
\def\Fr{\mbox{Fr}}
\def\gal{-Galois extension}
\def\hge{Hopf-Galois extension}
\def\ta{\tilde a}
\def\tb{\tilde b}
\def\tc{\tilde c}
\def\td{\tilde d}
\def\tx{\tilde x}
\def\ty{\tilde y}

\title{\vspace*{-15mm}{\large\bf EXPLICIT HOPF-GALOIS DESCRIPTION OF 
\boldmath $SL\sb{e\sp{\mbox{\tiny $\frac{2\pi i}{3}$}}}(2)$-INDUCED 
FROBENIUS HOMOMORPHISMS}\\
}
\author{\normalsize\sc
Ludwik D\c{a}browski\thanks{\small dabrow@sissa.it},~
Piotr M.~Hajac\thanks{\small On~leave from:
Department of Mathematical Methods in Physics, 
Warsaw University, ul.~Ho\.{z}a 74, Warsaw, \mbox{00--682~Poland}.
http://info.fuw.edu.pl/KMMF/ludzie\underline{~~}ang.html 
(E-mail: pmh@fuw.edu.pl)},~
Pasquale Siniscalco\thanks{\small sinis@sissa.it}
\\
\normalsize Scuola Internazionale di Studi Superiori Avanzati,\\
\normalsize Via Beirut 2-4, Trieste, \mbox{34014}~Italy.\\
\vspace*{-5mm}
}
\date{}
\maketitle

\begin{abstract}
The exact sequence of ``coordinate-ring'' Hopf algebras
\mbox{$A(SL(2,\IC)\!)\!\stackrel{Fr}{\rightarrow}\! A(SL_q(2)\!)
\!\rightarrow\! A(\! F)$} 
determined by the Frobenius map $Fr$, and the same way obtained 
exact sequence of (quantum) Borel subgroups,
are studied when $q$ is a cubic root of unity. 
An \slc-linear splitting of \slq\ making \slc\ a direct summand of \slq\
is constructed and used to prove that $A(SL_q(2))$
 is a faithfully flat $A(F)$-Galois extension of $A(SL(2,\IC))$.
A cocycle and coaction determining the
bicrossed-product structure of the upper-triangular (Borel)
 quantum subgroup of \slq\ are computed explicitly.
\end{abstract}

\section*{Introduction}

This work was inspired by the final remark in \cite{c-a} pointing to the
possible importance of quantum-group Frobenius homomorphisms in understanding
the (quantum) symmetry of the Standard Model. We focus our attention on the
cubic root of unity because it is the simplest non-trivial odd case and
 because, as advocated by A.~Connes, it might be the ``cubic symmetry"
that is to succeed the supersymmetry in physics.

In the present
 study of short exact sequences of quantum groups we adopt the 
functions-on-group point of view, which is dual to the 
universal-enveloping-algebra approach (see Paragraph~8.17 in~\cite{l}).
It is known \cite{a,ms} that Frobenius mappings at primitive odd roots of
unity allow us to view \slq\ as a faithfully flat \hge\ of \slc.
The main contribution of this paper is a construction of an \slc-linear splitting
of \slq\ making \slc\ a direct summand of \slq, and the computation of 
the cocycle-bicrossed-product structure of the analogous quantum extension of the
upper-triangular (Borel) subgroup of $SL(2,\IC)$. 
With the aim of attracting a diverse readership,
we write this article in a relatively self-contained down-to-earth
manner. We hope that, by exemplifying
certain concepts in a very tangible way, this note can serve as an
invitation to further study.

In the next two sections, we establish the basic language of this work
and review appropriate modifications of known
general results that we apply later to compute examples.

In Section~3, we reduce the task of computing the $A(F)$-coinvariants of \slq\
to finding a certain \slc-homomorphism. Just as Hopf-Galois extensions
generalise to a great extent the concept of a principal bundle, this homomorphism
generalises the notion of a section of a bundle. Thus we
derive an alternative proof that \slq\ is a faithfully flat \hge\ of \slc.

Section~4 and Section~5 are devoted to the study of the same kind Frobenius 
homomorphisms in the
Borel and Cartan cases. As the Hopf algebra $P_+:=\slq/\lan c\ran$ is pointed, we
can conclude that $P_+$ is a cleft \hge\ of $B_+:=\slc/\lan\bar{c}\ran$. 
We construct a family of cleaving maps 
\mbox{$\af/\lan\tc\ran=:H_+\st{\Phi_\nu}{\ra}P_+$}, calculate
an associated cocycle and weak coaction, and prove that $P_+$ has a 
non-trivial bicrossed-product structure. Our construction works for any
primitive odd root of unity.
The Cartan case (the off-diagonal generators put to zero) is commutative
and follows closely the Borel case pattern.

For the sake of completeness, in the final two sections
 we determine the integrals in and on \af,
prove the non-existence of the Haar measure on $F$, and show that the natural
\af-coinvariants of the polynomial algebra of the quantum plane at the cubic
root of unity form an algebra isomorphic with the algebra of polynomials on~$\IC^2$.
We also present corepresentations of~\af.

Throughout this paper we use Sweedler's notation (with the summation symbol
suppressed) for the coproduct ($\hD\ h=h\1\ot h\2$) and right coaction
($\dr\ p=p\0\ot p\1$). The unadorned tensor product stands for the tensor
product over a field $k$. (In the examples $k=\IC$.) 
The counit and antipode are denoted by \he\ and $S$ 
respectively, and $m$ is used to signify the multiplication in an algebra.
By the convolution product of two linear maps we understand 
$f*g:=m\ci(f\ot g)\ci\hD$, $(f*g)(h)=f(h\1)g(h\2)$. The convolution inverse
of $f$ is denoted by $f^{-1}$ and defined by $f*f^{-1}=\he=f^{-1}*f$.
We use $\hd_{kl}$ to denote the Kronecker delta.

\parskip0pt
\parindent0pt
\section{Preliminaries}

We begin by recalling basic definitions.

\bde\Label{hgdef} 
Let $H$ be a Hopf algebra, $P$ be a right $H$-comodule algebra, and 
$B:=P\sp{co H}:=\{p\in P\,|\; \dr\ p=p\ot 1\}$ the subalgebra of
right coinvariants. We say that
P is a (right) {\em \hge} (or $H$-Galois extension)
of $B$ iff the canonical left $P$-module
right $H$-comodule map 
$
can:=(m\sb P\ot id)\circ (id\ot\sb B \dr )\, :\; P\ot\sb B P\lra P\ot H
$
is bijective.
\ede

In what follows, we will use only right \hge s, and skip writing ``right"
for brevity.

\bde\Label{ffd} We say that $P$ is a faithfully flat $H$\gal\ of $B$
iff $P$ is faithfully flat as a right and left $B$-module.
(For a comprehensive review of the concept of faithful flatness see~\cite{b-n}.)
\ede

\bde\Label{cleft} 
An $H$\gal\ is called {\em cleft} iff there exists a  convolution
invertible linear map $\Phi : H\ra P$ satisfying 
$\dr\circ\Phi=(\Phi\ot id)\circ\hD$. We call $\Phi$ a cleaving map of $P$.
\ede

Note that, in general, $\Phi$ is {\em not} uniquely determined by its defining 
conditions. Observe also that a cleaving map
can always be normalised to be unital. Indeed, let $\tilde{\Phi}$ be  a cleaving map,
and $\tilde{\Phi}(1):=b$. By the colinearity, we have that $b\in B$, and
 the convolution invertibility entails that $b$ is invertible. Also, 
$b^{-1}\ot 1=b^{-1}\dr(bb^{-1})=b^{-1}b\dr(b^{-1})=\dr(b^{-1})$. 
It is straightforward to check that $\Phi:=b^{-1}\tilde{\Phi}$ is right colinear,
 convolution invertible and unital. Hence, without the loss of generality,
 we  assume $\Phi$ to be unital for the rest of this paper.
Let us also remark that a cleaving map
is necessarily injective:
\[
\llp m_P\ci(m_P\ot id)\ci(id\ot\Phi^{-1}\ot id)\ci(id\ot\hD)\ci\dr\ci\Phi\lrp(h)
=\Phi(h\1)\Phi^{-1}(h\2)h\3=h,~~~
\forall\, h\in H.
\]
\bde[\cite{pw}]\Label{es} 
A sequence of Hopf algebras (and Hopf algebra maps) 
\mbox{$B\st{j}{\ra} P\st{\pi}{\ra} H$} is called {\em exact} 
iff $j$ is injective and $\pi$ is the canonical surjection on
$H=P\slash Pj(B\sp+)P$, where $B\sp+$ denotes the augmentation ideal of $B$
(kernel of the counit map). 
\ede

When no confusion arises regarding the considered class of ``functions"
on quantum groups, one can use the above definition to define exact 
sequences
of quantum groups (see (1.6a) in \cite{pw}). 
In particular, we can view $F$ as a finite quantum group.
Further sophistication of 
the concept of a short exact sequence of quantum groups comes with the
following definition (cf.~\cite[p.23]{ad}):

\bde[p.3338 in \cite{s2}]\Label{ses} 
An exact sequence of Hopf algebras \mbox{$B\st{j}{\ra} P\st{\pi}{\ra} H$}
is called {\em strictly exact} iff $P$ is right faithfully flat over $j(B)$,
and $j(B)$ is a normal Hopf subalgebra of $P$, i.e., 
$\llp p\1 j(B)S(p\2)\cup S(p\1) j(B)p\2\lrp\inc j(B)$ for any $p\in P$.
\ede
See~\cite[Section~5]{m-a}
for short exact sequences of finite dimensional Hopf algebras.
\bre\em
Exact sequences of Hopf algebras should not be confused with exact sequences
of vector spaces: The exact sequence of groups 
\mbox{$\IZ\sb3\ra\IZ\sb6\ra\IZ\sb6\slash\IZ\sb3\cong\IZ\sb2$}
yields (by duality) an exact sequence of Hopf algebras which is {\em not}
an exact sequence in the category of vector spaces (or algebras).
\ere

Let us now provide a modification of Remark~1.2(1) in~\cite{s1} that allows us
to avoid directly verifying the faithful flatness condition in the 
proof of Proposition~\ref{coin}. We replace the faithful flatness assumption
by assuming the existence of a certain homomorphism. Its existence in the case
described in Proposition~\ref{coin} is proved through a calculation 
(Lemma~\ref{fr-linear}).

\ble\Label{coi}
 Let $P$ be a right $H$-comodule algebra and $C$ a {\em subalgebra}
of $P\sp{co H}$ such that the map
$
\psi : P\ot\sb CP\ni p\ot\sb Cp'\mapsto pp'\0\ot p'\1\in P\ot H
$
is bijective, and such that there exists a unital right $C$-linear homomorphism
$s:P\ra C$ (cf.~Definition~A.4 in~\cite{h}). 
Then $C=P\sp{co H}$, and $P$ is an $H$\gal\ of $C$.
\ele
\bpf
Note first that the map $\psi$ is well defined due to the assumption
$C\inc P\sp{co H}$. Now, let $x$ be an arbitrary element of $P\sp{co H}$.
Then
\beq\Label{dx}
1\ot\sb C x=\psi\sp{-1}(\psi(1\ot\sb C x))=\psi\sp{-1}(x\ot 1)=x\ot\sb C1\, .
\eeq
On the other hand, we know from Proposition~2.5 of~\cite{cq} that 
$P\ot\sb C(P\slash C)$ is isomorphic with 
$\mbox{Ker}(m\sb p:P\ot\sb CP\ra P)$. In particular, this isomorphism sends
$1\ot\sb C x -x\ot\sb C 1$ to $1\ot\sb C [x]\sb C\in P\ot\sb C (P\slash C)$.
Remembering (\ref{dx}) and applying first $s\ot\sb C id$ and then the 
multiplication map to 
$1\ot\sb C [x]\sb C$, we obtain $[x]\sb C=0$, i.e.\ $x\in C$, as needed.
\epf

\bre\Label{flat}\em 
Observe that the assumption of the existence of a unital right $C$-linear 
homomorphism $s:P\ra C$ can be replaced by the assumption that $P\slash C$ is
flat as a left $C$-module. Indeed, we could then view $C\ot\sb C (P\slash C)$
as a submodule of $P\ot\sb C (P\slash C)$, and consequently $1\ot\sb C [x]\sb C$
as an element of the former. Now one could directly apply the multiplication
map to $1\ot\sb C [x]\sb C$ and conclude the proof as before.
\ere

It is well known that cleft \hge s and crossed products are equivalent notions.
Once we have a cleaving map $\Phi$,  we can determine
the cocycle and cocycle action that define the crossed product structure (see
\cite[Section~4]{bcm}, \cite[Definition~6.3.1]{m-s}) from the 
following formulas respectively \cite[p.273]{s3}:  
\beq\Label{act}
h\triangleright\sb\Phi b:=\Phi(h\1)b\Phi\sp{-1}(h\2)\in P^{co H}
\eeq\beq\Label{coc}
\sigma\sb\Phi(h\ot l):=\Phi(h\1)\Phi(l\1)\Phi\sp{-1}(h\2 l\2)\in P^{co H}\; ,
\eeq
where $h,l\in H,\; b\in P\sp{co H}$. On the other hand, 
with the help of $\Phi$
we can construct a unital left $B$-module homomorphism 
$s\sb\Phi:P\ra B$ by the formula 
\beq\Label{sphi}
s\sb\Phi:=m\circ (id\ot\Phi\sp{-1})\circ\dr\, .
\eeq
The homomorphism $s\sb\Phi$ generalises the notion of a section of a principal
bundle just as $\Phi$ generalises the concept of a trivialisation of a principal
bundle (see the end of Section~4 here and Remark~2.6 in~\cite{h}). The following
straightforward-to-prove lemma allows one to compute $\sigma\sb\Phi$ by taking
advantage of $s\sb\Phi$. It seems to be a more convenient way of calculating
$\sigma\sb\Phi$ whenever $\hD\ot\hD$ is more complicated than \dr.  We
will use it to compute a cocycle of the cleft extension describing an exact
sequence of (quantum) Borel subgroups.

\ble[cf.~Lemma 2.5 in \cite{h}]\Label{cycle} 
Let $P$ be a cleft $H$\gal\ of $B$ and $\Phi$ a
cleaving map. Then $\sigma\sb\Phi=s\sb\Phi\circ m\circ (\Phi\ot\Phi)\, .$
\ele

Finally, let us observe that with the help of the translation map 
(e.g., see~\cite{b-t}) 
\[
\tau:H\ra P\ot_BP,\;
\tau(h):=can^{-1}(1\ot h)=:h^{(1)}\ot_Bh^{(2)}
\]
 (summation suppressed),
we can solve formula (\ref{sphi}) for $\Phi$. Indeed, 
\[
(id*_\tau s_\Phi)(h):=h^{(1)}s_\Phi(h^{(2)})
=h^{(1)}{h^{(2)}}_{(0)}\Phi^{-1}({h^{(2)}}_{(1)})
=\llp m\ci(id\ot\Phi^{-1})\ci can\lrp(h^{(1)}\ot_Bh^{(2)})=\Phi^{-1}(h),
\]
whence $\Phi=(id*_\tau s_\Phi)^{-1}$. \\

\section{Principal homogenous extensions}

Let $P$ be a Hopf algebra and a $(P/I)$\gal\ of $B$ for
the coaction 
\[
\dr:=(id\ot\pi)\ci\hD,\;\;\; P\st{\pi}{\ra}P/I,
\]
 where $I$ is
a Hopf ideal of~$P$. Then we call $P$ is a {\em principal homogenous extension} 
of~$B$. First we recall a theorem\footnote{
 We owe it   to Peter Schauenburg.
}
 which shows the structure of the Hopf ideal~$I$. 

\bth[cf.~Lemma~5.2 in \cite{bm}]\Label{b+p}
Let $P$ be a Hopf algebra and $I$ a Hopf ideal of $P$. Then 
$P$ is a $(P/I)$-principal homogenous extension of $B$
{\em if and only if} $I=B^+P$, where $B:=P^{co(P/I)},\; 
B^+:=B\cap\mbox{\em Ker}\,\he$.
\ethe
{\it Proof.}~
Assume first that $I=B^+P$. Taking advantage of (\ref{bp}), for any
$b\in B^+$, $p\in P$, we have:
\beq\Label{sbp}
S(b\1 p\1)\ot_Bb\2 p\2=S(b\1 p\1)b\2\ot_Bp\2=S(p\1)\he(b)\ot_Bp\2=0.
\eeq
Hence we have a well-defined map $\wp:P\ot(P/I)\ra P\ot_BP,\;
\wp(p\ot [p']_I):=pS(p'\1)\ot_Bp'\2$. It is straightforward to verify
that $\wp$ is the inverse of the canonical map $can$. Consequently,
$P$ is a \hge.\\

To show the converse, let us first prove the following:

\ble\Label{homo}
Let $P$, $I$ and $B$ be as above. Then $B\inc P$ is a $(P/I)$\gal\
{\em if and only if}
$\llp\pi\sb B\circ(S\ot id)\circ\hD\lrp(I)=0$, where 
$\pi\sb B:P\ot P\ra P\ot\sb B P$ is the canonical surjection.
\ele
\bpf
If $P$ is a $(P/I)$\gal\ of $B$, then 
we have the following short exact sequence (see the proof of Proposition~1.6
in~\cite{h}):
\beq\Label{seq}
0\ra P(\hO\sp1\! B)P\hookrightarrow 
P\ot P\st{T\sb R}{\ra}P\ot P\slash I\ra 0.
\eeq
Here $\hO\sp1\! B:=\mbox{Ker}\,(m:B\ot B\ra B)$ and 
$T\sb R=(m\ot\pi)\circ(id\ot\hD)$. One can check
that $\llp T\sb R\circ(S\ot id)\circ\hD\lrp(I)=0$. Hence, it follows from 
the exactness of (\ref{seq}) that 
$\llp (S\ot id)\circ\hD\lrp(I)\inc P(\hO\sp1\! B)P$. Consequently,
$\llp \pi\sb B\circ(S\ot id)\circ\hD\lrp(I)=0$ due to the exactness of 
the sequence
\[
0\ra P(\hO\sp1\! B)P\hookrightarrow 
P\ot P\st{\pi\sb B}{\ra}P\ot\sb B P\ra 0\, .
\]
To prove the converse, one can proceed as in the considerations 
preceding this lemma.
\epf

\bco\Label{bpco}
Let $B\inc P$ be a $(P/I)$\gal\ as above. Then the translation map
is given by the formula: $\tau([p]_I):=S(p\1)\ot_Bp\2\,$.
\eco

Assume now that $P$ is a $(P/I)$\gal\ of $B$. It follows from the above
corollary and (\ref{sbp}) that $\tau([B^+P]_I)=0$. 
Hence, by the injectivity
of $\tau$, we have $B^+P\inc I$. Furthermore, we have a well-defined
map $can':P\ot_BP\ra P\ot(P/B^+P),\; p\ot_Bp'\to pp'\1\ot [p'\2]_{B^+P}\,$.
Indeed, taking again  advantage of (\ref{bp}), we obtain
\begin{eqnarray*}
p\ot bp'&\to& pb\1 p'\1\ot [(b\2-\he(b\2))p'\2+\he(b\2)p'\2]_{B^+P}\\
&=& pb\1 p'\1\ot\he(b\2)[p'\2]_{B^+P}\\
&=& pbp'\1\ot [p'\2]_{B^+P}\, ,
\end{eqnarray*}
and $pb\ot p'\to pbp'\1\ot [p'\2]_{B^+P}\,$. Reasoning as in the first part
of the proof, we can conclude that $can'$ is bijective. We have the following
commutative diagram:

\[
\begin{array}{ll}
P\ot_BP  & \mathop{-\hspace{-6pt}\longrightarrow}\limits^{can'} P\ot(P/B^+P)  \\

~\, ^{id} \Big\downarrow  &  \,~~~~~~~~~~\Big\downarrow\ ^{id\ot\ell} \\ 

P\ot_BP & \mathop{-\hspace{-6pt}\longrightarrow}\limits^{can} P\ot(P/I)\, , \\
\end{array}
\]\ \\
where $\ell([p]_{B^+P}):=[p]_{I}$. (Recall that we have already showed that
$B^+P\inc I$, so that $\ell$ is well defined.) 
It follows from the commutativity of the diagram that $id\ot\ell$
is bijective. In particular, we have that $\ell$ is injective, and therefore
$I\inc B^+P$, as needed. \hfill{$\rule{7pt}{7pt}$}\\

Let us now prove the following left-sided
version of a result by Y.Doi and A.Masuoka (see \cite{md92} or 
\cite[Proposition~3.8]{m-a}):

\bth\Label{cinv}
Let $P$ be a $(P/I)$-principal homogenous extension of $B$.
 Then $P$ is cleft
{\em if and only if} there exists a convolution invertible  left
$B$-module homomorphism $\Psi:P\ra B$.
\ethe\bpf
Assume first that $P$ is cleft. Let $\Phi$ be a cleaving map. Then  
$\Psi:=s_\Phi$ (see~(\ref{sphi})) is left $B$-linear. 
 Moreover, it can be directly verified that 
$\Psi^{-1}:P\ra B,\;\Psi^{-1}(p):=\Phi(\pi(p\1))S(p\2)$, 
(see~\cite[Definition 3.2.13(3)]{ad}) is the convolution inverse of $\Psi$.\\

Conversely, assume that we have $\Psi:P\ra B$ with the required properties.
To define $\Phi$ in terms of $\Psi$, first we need to derive certain property
of $\Psi^{-1}$.

\ble\Label{psi-}
Let $\Psi:P\ra B$ be a homomorphism as described in Theorem~\ref{cinv}. Then
\linebreak
$\Psi^{-1}(b\1 p)b\2=\he(b)\Psi^{-1}(p)$ for any $b\in B,\; p\in P$. 
\ele\bpf
Note first that $b\in B$ implies $b\1\ot b\2\in P\ot B$. Indeed,
\bea\Label{bp}
(id\ot\dr)(b\1\ot b\2)
&=&\llp(id\ot id\ot\pi)\ci(id\ot\hD)\ci\hD\lrp(b)\nonumber\\
&=&\llp(id\ot id\ot\pi)\ci(\hD\ot id)\ci\hD\lrp(b)\nonumber\\
&=&\llp(\hD\ot id)\ci\dr\lrp(b)\nonumber\\
&=&b\1\ot b\2\ot 1.
\eea
Taking advantage of this fact, for any $b\in B,\; p\in P$, we obtain
\begin{eqnarray*}
\he(b)\Psi^{-1}(p)
&=&\Psi^{-1}(b\1 p\1)\Psi(b\2 p\2)\Psi^{-1}(p\3)\\
&=&\Psi^{-1}(b\1 p\1)b\2\Psi(p\2)\Psi^{-1}(p\3)\\
&=&\Psi^{-1}(b\1 p)b\2\, ,
\end{eqnarray*}
as claimed.
\epf\ \\

On the other hand, we know (see Theorem~\ref{b+p}) that, since $P$ is
a $(P/I)$\gal, $I=B^+P$. Furthermore, with the help of Lemma~\ref{psi-},
we can directly show that $(\Psi^{-1}*id)(B^+P)=0$. Hence we have a 
well-defined map $\Phi:P/I\ra P,\; \Phi(\pi(p)):=(\Psi^{-1}*id)(p)=
\Psi^{-1}(p\1)p\2\,$.  We also have:
\[
(\dr\ci\Phi\ci\pi)(p)
=(\Psi^{-1}(p\1)\ot 1)(id\ot\pi)(\hD p\2)
=\Psi^{-1}(p\1)p\2\ot\pi(p\3)
=\Phi(\pi(p\1))\ot\pi(p\2)\, ,
\]
i.e., $\Phi$ is colinear. As expected from the general discussion in the
previous section, the formula for the convolution inverse of $\Phi$ is
$\Phi^{-1}=id*_\tau\Psi$. In our case we know that the formula for
the translation map is $\tau(\pi(p))=S(p\1)\ot_Bp\2$ (see~\cite[p.294]{s1}
and Corollary~\ref{bpco}). Thus we obtain: $\Phi^{-1}(\pi(p))=S(p\1)\Psi(p\2)\,$.
It can be directly checked that $\Phi^{-1}$ is indeed the convolution
inverse of $\Phi$. \hfill{$\rule{7pt}{7pt}$}\\

In the spirit of~\cite[p.47]{ad},
we say that $P$ is {\em cocleft} iff there exists a  convolution
invertible left $B$-module map (retraction) $\Psi:P\ra B$. 
(See \cite[Definition~2.2]{md92} for the right-sided version.) We call $\Psi$
a cocleaving map.

\ble\Label{norm}
Let $P$ be $P/I$-\phe\ of~$B$. Any cleaving and any cocleaving map of such an
extension can always be normalised to be both unital and counital.
\ele
\bpf
We already know from the previous section that a cleaving map can always be made
unital. Similarly, for any cocleaving map $\widetilde{\Psi}:P\ra B$, the map defined
by $\check{\Psi}(p)=\widetilde{\Psi}(p)\widetilde{\Psi}(1)^{-1}$ 
is a unital cocleaving
map. Here the invertibility of $\widetilde{\Psi}(1)$ follows from the convolution 
invertibility of~$\widetilde{\Psi}$, and 
$\check{\Psi}^{-1}(p)=\widetilde{\Psi}(1)\widetilde{\Psi}^{-1}(p)$ is the formula
for the convolution inverse of~$\check{\Psi}$. We can describe $\check{\Psi}$ as
the composite mapping:
\[
P
~\st{id\ot 1}{-\!\!\!-\!\!\!-\!\!\!-\!\!\!\lra}~
P\ot k
~\st{\widetilde{\Psi}\ot\widetilde{\Psi}(.)^{-1}}{-\!\!\!-\!\!\!-\!\!\!-\!\!\!\lra}~
B\ot B
~\st{m_B}{-\!\!\!-\!\!\!-\!\!\!-\!\!\!\lra}~
B.
\]
Formally ``dualising" this sequence and exchanging factors in the tensor product
one obtains:
\[
H
~\st{\Delta}{-\!\!\!-\!\!\!-\!\!\!-\!\!\!\lra}~
H\ot H
~\st{(\epsilon\ci f^{-1})\ot f}{-\!\!\!-\!\!\!-\!\!\!-\!\!\!\lra}~
k\ot P
~\st{m}{-\!\!\!-\!\!\!-\!\!\!-\!\!\!\lra}~
P.
\]
This suggests that one can counitalise a unital cleaving map 
$\check{\Phi}:P/I\ra P$ by the formula\footnote{
We owe this formula to a discussion with Tomasz Brzezi\'nski.
}
$\Phi(h)=\he(\check{\Phi}^{-1}(h\1))\check{\Phi}(h\2)$.
Indeed, it is straightforward to check that the thus defined map is counital,
unital, colinear and convolution invertible with the convolution inverse
given by $\Phi^{-1}(h)=\check{\Phi}^{-1}(h\1)\he(\check{\Phi}(h\2))$.
To complete the proof it suffices to check that if $\check{\Psi}:P\ra B$
is a unital cocleaving map, then the map defined by the formula
$\Psi(p)=\check{\Psi}(p\1)\he(\check{\Psi}^{-1}(p\2))$ is unital, counital
and cocleaving. The first two properties are immediate. Furthermore, it is
straightforward to verify that 
$\Psi^{-1}(p)=\he(\check{\Psi}(p\1))\check{\Psi}^{-1}(p\2)$ defines the
convolution inverse of~$\Psi$. It remains to make sure that $\Psi$ is left
$B$-linear. To this end, taking advantage of the fact that
$b\in B\Rightarrow b\1\ot b\2\in P\ot B$ (see~(\ref{bp})) and Lemma~\ref{psi-},
we compute:
\begin{eqnarray*}
\Psi(bp)\!\!\!\!\!\!\!
&& =
\check{\Psi}(b\1 p\1)\he(\check{\Psi}^{-1}(b\2 p\2))
\\ && = 
\check{\Psi}(b\1 p\1)\he(\check{\Psi}^{-1}(b\2 p\2)b\3)
\\ && = 
\check{\Psi}(b\1 p\1)\he(\he(b\2)\check{\Psi}^{-1}(p\2))
\\ && = 
\check{\Psi}(bp\1)\he(\check{\Psi}^{-1}(p\2))
\\ && = 
b\Psi(p).
\end{eqnarray*}
In the last step we used the assumption that $\check{\Psi}$ is left $B$-linear. 
\epf\\

\bco\Label{psico}
Let $P$ be $P/I$-\phe\ of~$B$. Then the following statements are equivalent:\\

1. $P$ is cleft.

2. $P$ is cocleft.

3. There exists a unital and counital convolution invertible right
$P/I$-colinear map \mbox{$\Phi: P/I\ra P$}.

4. There exists a unital and counital convolution invertible left
$B$-linear map $\Psi: P\ra B$.\\
 
Also, we have a one-to-one correspondence between
the unital-counital cleaving  and the unital-counital
cocleaving maps of $P$. The formula
\[
\Phi\longmapsto\Psi:=s_\Phi:=m\ci(id\ot\Phi^{-1})\ci\dr
\]
defines the desired bijection. Its inverse is given by
\[
\Psi\longmapsto\Phi,\;\;\;\Phi(\pi(p)):=(\Psi^{-1}*id)(p)=\Psi^{-1}(p\1)p\2\, .
\]
\eco

Observe that our considerations are very similar to those on p.47 and p.50
in~\cite{ad}. Here, however, we do not assume that the algebra of coinvariants
is a Hopf algebra.\\

\section{\boldmath $A(SL\sb{e\sp{\mbox{\tiny $\frac{2\pi i}{3}$}}}(2))$
as a faithfully flat Hopf-Galois extension}

Recall that $A(SL_q(2))$ is a complex Hopf algebra
generated by $1,\, a,\, b,\, c,\, d$, satisfying the 
following relations: 
\begin{eqnarray*}
&& ab=q ba~,~~ ac=qca~,~~ bd=qdb~,~~ bc = cb~,~~  cd =qdc~, \\ 
&& ad-da=(q-q^{-1})bc~, ~~ ad-qbc=da-q^{-1 }bc=1\, ,  
\end{eqnarray*}
where $q\in\IC\setminus\{0\}$. 
The comultiplication $\Delta$, counit $\varepsilon$, and antipode
$S$ of \slq\ are defined by the following formulas:
\[ 
\Delta {\scriptsize \left( \ba{cc} a & b \\ c & d \ea \right)} =
{\scriptsize \left( \ba{cc} a\ot 1 & b\ot 1 \\ c\ot 1 & d\ot 1 \ea \right)}  
{\scriptsize \left( \ba{cc} 1\ot a & 1\ot b \\ 1\ot c & 1\ot d \ea \right)}\, ,\;\;
\varepsilon {\scriptsize \left( \ba{cc} a & b \\ c & d \ea \right)}=
{\scriptsize \left( \ba{cc} 1 & 0 \\ 0 & 1 \ea \right)}\, ,\;\;
S{\scriptsize \left( \ba{cc} a & b \\ c & d \ea \right)}=
{\scriptsize \left( \ba{cc} d &-q^{-1} b \\ -qc & a \ea \right)}\, .
\]

Let us now establish some notation (e.g., see Section~IV.2 in~\cite{k}): 
\bea
&& 
(k)_q := 1+q+ \dots + q^{k-1}= \frac{ q^{k}-1}{q-1}
\, ,\;~ k\in {\IZ},\; k>0\, ;
\nonumber\\ &&
(k)_q ! := (1)_q(2)_q\dots (k)_q = \frac{ (q-1)(q^2-1) \dots 
 (q^k-1)}{(q-1)^k}\, ,\;~ (0)_q!:= 1\, ;
\nonumber\\ &&
\left({\scriptsize 
\ba{c}\!\!\! k\!\!\! \\ \!\!\! i\!\!\! \ea }\right)_q :=  
 \frac{ (k)_q !}{(k -i)_q ! (i)_q !}\, ,\;~ 0\leq i \leq k\, .\nonumber
\eea
The above defined $q$-binomial coefficients satisfy the following equality:
\[
(u+v)^k=\sum_{l=0}^{k}\left({\scriptsize 
\ba{c}\!\!\! k\!\!\! \\ \!\!\! l\!\!\! \ea }\right)_q u^lv^{k-l},
\]
where $uv=q^{-1}vu$. Now, if 
$(T_{ij})={\scriptsize\pmatrix{a&b\cr c&d\cr}}$, then
\[
\hD T_{ij}^k=(T_{i1}\te T_{1j}+T_{i2}\te T_{2j})^k
=\sum_{l=0}^{k}\left({\scriptsize 
\ba{c}\!\!\! k\!\!\! \\ \!\!\! l\!\!\! \ea }\right)_{q^{-2}}
T_{i1}^lT_{i2}^{k-l}\ot T_{1j}^lT_{2j}^{k-l}\, .
\]
{\em For the rest of this paper we put $q=e\sp{\frac{2\pi i}{3}}$.} Obviously,
we now have  $q^{-2}=q$, and the comultiplication on the basis elements
of \slq\ (see Lemma~1.4 in~\cite{mmnnu}, Exercise~7 on p.90 in~\cite{k})
 is given by:
\bea \Label{co}
\Delta(a^pb^rc^s) = \sum_{\lambda,\mu,\nu=0}^{p,r,s} 
\left({\scriptsize \ba{c}\!\!\! p\!\!\! \\ \!\!\! \lambda\!\!\! \ea }\right)_q
\left({\scriptsize \ba{c}\!\!\! r\!\!\! \\ \!\!\! \mu\!\!\! \ea }\right)_q
\left({\scriptsize \ba{c}\!\!\! s\!\!\! \\ \!\!\! \nu\!\!\! \ea }\right)_q
 a^{p-\lambda} b^{\lambda} a^{\mu} b^{r-\mu} c^{s-\nu} d^{\nu} 
\otimes a^{p- \lambda} c^{\lambda} b^{\mu} d^{r-\mu} a^{s-\nu} c^{\nu} 
\, ,\nonumber\\
\Delta(b^kc^ld^m) = \sum_{\lambda,\mu,\nu=0}^{k,l,m} 
\left({\scriptsize \ba{c}\!\!\! k\!\!\! \\ \!\!\! \lambda\!\!\! \ea }\right)_q
\left({\scriptsize \ba{c}\!\!\! l\!\!\! \\ \!\!\! \mu\!\!\! \ea }\right)_q
\left({\scriptsize \ba{c}\!\!\! m\!\!\! \\ \!\!\! \nu\!\!\! \ea }\right)_q
a^{\lambda} b^{k-\lambda} c^{l-\mu} d^{\mu} c^{\nu} d^{m-\nu} 
\otimes b^{\lambda} d^{k-\lambda} a^{l- \mu} c^{\mu} b^{\nu} d^{m-\nu}\, , 
\eea
where $m$ is a positive integer and $p,r,s,k,l$  are non-negative integers.\\

Following Chapter~7 of \cite{pw} and Section~4.5 of \cite{m-yu} (cf.~Section~5 in 
\cite{t} and the end of Part~I of~\cite{c-p}), we take the Frobenius mapping
\beq\Label{fr}
Fr:\slc\ni \bar{T}\sb{ij}\longmapsto T\sp3\sb{ij}\in\slq\, ,\;\; i,j\in\{1,2\},
\eeq
to construct the exact sequence of Hopf algebras
\beq\Label{frs}
\slc\st{Fr}{\lra}\slq\st{\pi\sb F}{\lra}\af\, .
\eeq
Here $\af=\slq\slash\lan T\sp3\sb{ij}-\hd\sb{ij}\ran\, ,\;\; i,j\in\{1,2\}$,
and $\pi\sb F$ is the canonical surjection. The following proposition
determines a basis of \af\ and shows that \af\ is 27-dimensional.

\bpr\Label{f-basis} 
Define $\ta:=\pi_F(a),\;\tb:=\pi_F(b),\;\tc:=\pi_F(c),\;\td:=\pi_F(d)$. 
Then the set $\{\tilde{a}^p \tilde{b}^r\tilde{c}^s\}\sb{p,r,s\in\{0,1,2\}}$ 
is a basis of \af .
\epr\bpf
Since $\td=\ta^2(1+q\tb\tc)$, the monomials 
$\tilde{a}^p \tilde{b}^r\tilde{c}^s,\; p,r,s\in\{0,1,2\}$, 
span \af. 
Guided by the left action of \af\ on itself,
we define a 27-dimensional representation 
 $\varrho : \af\ra \mbox{End}(\IC\sp3\ot\IC\sp3\ot\IC\sp3)$ 
 by the following formulas:
\bea 
\varrho(\tilde{a}) & = & {\bf J} \ot {\bf I\sb3} \ot {\bf I\sb3}~,
\nonumber \\
      \varrho(\tilde{b}) & = & {\bf Q} \ot {\bf N} \ot {\bf I\sb3}~, 
\nonumber \\
\varrho(\tilde{c}) & = & {\bf Q} \ot {\bf I\sb3} \ot {\bf N}~, \nonumber
\eea
 where 
\[
\bf{J} =  \left( \ba{ccc} 0 & 0 & 1 \\ 1 & 0 & 0 \\ 0  & 1 & 0 \ea 
              \right),~
\bf{Q} = \left( \ba{ccc} 1 & 0 & 0 \\ 0 & q^{-1} & 0 \\ 0  & 0 & q^{-2} \ea 
   \right),~
\bf{N} =\left( \ba{ccc} 0 & 0 & 0 \\ 1 & 0 & 0 \\ 0  & 1 & 0 \ea 
\right),~
\bf{I\sb3} =\left( \ba{ccc} 1 & 0 & 0 \\ 0 & 1 & 0 \\ 0  & 0 & 1 \ea 
\right). 
\]
It is straightforward to check that $\varrho$ is well defined.
Assume now that
$
\sum\sb{p,r,s\in\{0,1,2\}} \al_{prs} \tilde{a}^p\tilde{b}^r\tilde{c}^s=0. 
$ 
Applying $\varrho$, we obtain 
\beq\Label{jqn} 
\sum\sb{p,r,s\in\{0,1,2\}}\al_{prs} {\bf J}^p{\bf Q}^{r+s} \ot {\bf N}^r 
\ot{\bf N}^s=0\, .
\eeq
On the other hand, let us consider the linear functionals 
\[ 
h^{klm}: M_3(\IC)^{\ot^3}\ra \IC\, ,\;\;\;  
h^{klm} (A \ot B \ot C) := A_{k0}B_{l0}C_{m0}\, ,
\;\;\; k,l,m\in\{0,1,2\}\, , 
\] 
 where we number the rows and columns of matrices by 0,1,2.
From (\ref{jqn}) we can conclude that 
 \[ 
h^{klm}\left(\sum\sb{p,r,s\in\{0,1,2\}}\al_{prs} {\bf J}^p{\bf Q}^{r+s} 
\ot {\bf N}^r
\ot {\bf  N}^s\right)= 0, ~~ \forall k,l,m\in\{0,1,2\} ~. 
\] 
Consequently, since 
$ 
h^{klm}( {\bf J}^p{\bf Q}^{r+s} \ot 
{\bf N}^r \ot {\bf N}^s)=  \delta_{pk}\delta_{rl}\delta_{ms}\, 
$, 
we have that $\al_{prs} = 0$, for any $p,r,s$. Hence 
$\tilde{a}^p\tilde{b}^r\tilde{c}^s$ are linearly independent, as claimed.
\epf

\bco[cf.~Section~3 in~\cite{s-a}] The representation 
$\varrho : \af\ra \mbox{\em End}(\IC\sp3\ot\IC\sp3\ot\IC\sp3)$
defined above is faithful.
\eco

\bre\em
Observe that we could equally well consider a representation with
 ${\bf Q}$ replaced by ${\bf Q}^{-1}$, ${\bf J}$ by ${\bf J}^{t}$ and
${\bf N}$ by ${\bf N}^{t}$, where $\phantom{M}\sp t$ 
denotes the matrix transpose.
\ere

With the help of duality between functions-on-group and universal-enveloping
algebra pictures, it can be shown (\cite[Proposition~3.4.5]{a},~\cite{ms})
that $\fr=\slq\sp{A(F)}$. This can also be concluded from the fact that
\slq\ is Noetherian (see \cite[Theorem~IV.4.1, Proposition~I.8.2]{k}) 
and the combination 
of \cite[Theorem~3.3]{s2}, \cite[Remark~1.2(1)]{s1} and  
\cite[Remark~1.6(1)]{s2}. 
Note that Theorem~3.3 in~\cite{s2} establishes the faithful flatness
of \slq\ over \slc. After identifying \slc\ with $\slq\sp{co A(F)}$, we
can use Theorem~1.3 in~\cite{s3} (see~\cite{kt}) to infer
that \slq\ is finitely generated projective over \slc\ 
(cf.~\cite[Proposition~1.7]{dl}).
In what follows, we provide a direct proof which does not invoke the duality. 

\bpr\Label{coin}
The algebra $A(SL(2,\IC ))$ of polynomial functions on $SL(2,\IC )$ is 
isomorphic
(via the Frobenius map) to the subalgebra $A(SL\sb q(2))^{co A(F)}$ of all
right coinvariants. 
\epr
\bpf
The claim of the proposition follows immediately from the lemma below,
\cite[Lemma~1.3(1)]{s1} and Lemma~\ref{coi}. 
(From these lemmas one can also conclude that
\slq\ is an \af\gal\ of \fr .)

\ble \Label{fr-linear} 
Let $p,r,s,k,l,m\in\IN_0\, ,\; m>0$.
The linear map $s :\slq \ra \fr$ defined by the formulas
\bea
&&
s(a^pb^rc^s)=
\left\{\begin{array}{ll}
a^pb^rc^s & \mbox{\em when $p,r,s$ are divisible by 3}\\
0 & \mbox{\em otherwise}\, ,
\end{array}\right.
\nonumber\\ &&
s(b^kc^ld^m)=
\left\{\begin{array}{ll}
b^kc^ld^m & \mbox{\em when $k,l,m$ are divisible by 3}\\
0 & \mbox{\em otherwise}\, ,\nonumber
\end{array}\right.
\eea
is a unital \fr -homomorphism.
\ele
\bpf 
The unitality is obvious. Next, as \fr\ is a central subalgebra of \slq\
(see Theorem~5.1.(a) in~\cite{m-yu}),
 the left and right \fr-module structure of \slq\ coincide. Now,
we want to show that $s(f\omega)=fs(\omega)$, for any $f\in \fr$ and 
$\omega \in \slq$. In terms of the basis of \slq , we have a natural 
decomposition $f=f^1 +f^2$, $\omega=\omega^1+\omega^2$, where
$f^1=\sum f^1_{prs}a^{3p}b^{3r}c^{3s}$, 
$f^2=\sum\sb{m>0} f^2_{klm}b^{3k}c^{3l}d^{3m}$, 
$\omega^1=\sum \omega^1_{\alpha\beta\gamma}a^{\alpha}b^{\beta}c^{\gamma}$, 
$\omega^2=\sum\sb{\nu>0} \omega^2_{\lambda\mu\nu}b^{\lambda}c^{\mu}d^{\nu}$. 
Unless otherwise specified, we sum here over non-negative integers.
It is straightforward to see that $s(f^1\omega^1)=f^1s(\omega^1)$
 and $s(f^2\omega^2)=f^2s(\omega^2)$. We will demonstrate 
that $s(f^2\omega^1)=f^2s(\omega^1)$. 
We have: 

\bea
&& 
f^2\omega^1 =\sum\sb{m>0} f^2_{klm}b^{3k}c^{3l}d^{3m} \sum 
\omega^1_{\alpha\beta\gamma}a^{\alpha}b^{\beta}c^{\gamma} 
 = \sum\sb{m>0}  f^2_{klm} \omega^1_{\alpha\beta\gamma}
d^{3m}a^{\alpha}b^{3k+\beta}c^{3l+\gamma}\phantom{wwwwwwwwwww}
\nonumber \\ && \nonumber \\ && 
= \sum_{3m>\alpha}  f^2_{klm} \omega^1_{\alpha\beta\gamma}
d^{3m-\alpha}d^{\alpha}a^{\alpha}b^{3k+\beta}c^{3l+\gamma} + 
\sum_{0<3m\le \alpha}  f^2_{klm} \omega^1_{\alpha\beta\gamma}
d^{3m}a^{3m}a^{\alpha-3m}b^{3k+\beta}c^{3l+\gamma}
\nonumber \\ && \nonumber \\
&& = \sum_{3m>\alpha}  f^2_{klm} \omega^1_{\alpha\beta\gamma}
d^{3m-\alpha}p\sb\alpha(b,c)b^{3k+\beta}c^{3l+\gamma} + 
\sum_{0<3m\le \alpha}  f^2_{klm} \omega^1_{\alpha\beta\gamma}
a^{\alpha-3m}p_m(b^3,c^3)b^{3k+\beta}c^{3l+\gamma}.\nonumber  
\eea\

Here, due to the relation $da=1+q\sp{-1}bc$, the monomials
$d\sp\alpha a\sp\alpha =:p\sb\alpha(b,c)$ and $d^{3m}a\sp{3m}=:p_m(b^3,c^3)$
are polynomials in $b,c$ and $b^3,c^3$ respectively. Applying $s$ yields:

\bea
s(f^2\ho^1)&=&\sum_{m>\lambda}  f^2_{klm} \omega^1_{3\lambda,\beta,\gamma}
s\llp d^{3(m-\lambda)}p\sb{3\lambda}(b,c)b^{3k+\beta}c^{3l+\gamma}\lrp 
\nonumber \\ &+& 
\sum_{0<m\le\lambda}  f^2_{klm} \omega^1_{3\lambda,3\mu,3\nu}
s\llp a^{3(\lambda-m)}p_m(b^3,c^3)b^{3(k+\mu)}c^{3(l+\nu)}\lrp
\nonumber \\ &=&
\sum_{m>\lambda}  f^2_{klm} \omega^1_{3\lambda,\beta,\gamma}
s\llp d\sp{3\lambda}a\sp{3\lambda}b^{3k+\beta}c^{3l+\gamma}d^{3(m-\lambda)}\lrp 
\nonumber \\ &+& 
\sum_{0<m\le\lambda}  f^2_{klm} \omega^1_{3\lambda,3\mu,3\nu}
a^{3(\lambda-m)}p_m(b^3,c^3)b^{3(k+\mu)}c^{3(l+\nu)}
\nonumber \\ &=&
\sum_{m>\lambda}  f^2_{klm} \omega^1_{3\lambda,3\mu,3\nu}
d\sp{3\lambda}a\sp{3\lambda}b^{3(k+\mu)}c^{3(l+\nu)}d^{3(m-\lambda)}\lrp 
\nonumber \\ &+& 
\sum_{0<m\le\lambda}  f^2_{klm} \omega^1_{3\lambda,3\mu,3\nu}
a^{3(\lambda-m)}d^{3m}a^{3m}b^{3(k+\mu)}c^{3(l+\nu)}\, .\nonumber
\eea\ 

On the other hand, we have:
\bea 
 f^2s(\omega^1) &=& \sum\sb{m>0} f^2_{klm}b^{3k}c^{3l}d^{3m} \sum 
\omega^1_{3\lambda 3\mu 3\nu}a^{3\lambda}b^{3\mu}c^{3\nu} 
\nonumber \\ &=& 
\sum\sb{m>0} f^2_{klm} \omega^1_{3\lambda 3\mu 3\nu}
d^{3m}a^{3\lambda}b^{3k+3\mu}c^{3l+3\nu}
\nonumber \\ &=&
\sum_{m>\lambda} f^2_{klm} \omega^1_{3\lambda 3\mu 3\nu}
d^{3m-3\lambda}d^{3\lambda}a^{3\lambda}b^{3k+3\mu}c^{3l+3\nu} 
\nonumber \\ &+&
\sum_{0<m\le \lambda}  f^2_{klm} \omega^1_{3\lambda 3\mu 3\nu}
d^{3m}a^{3m}a^{3\lambda-3m}b^{3k+3\mu}c^{3l+3\nu}.\nonumber 
\eea\ 

Hence $s(f^2\omega^1)=f^2s(\omega^1)$, as needed. The remaining equality
$s(f^1\omega^2)=f^1s(\omega^2)$ can be proved in a similar manner.
\epf\\ \ 

Note that it follows from the above lemma that $P=B\oplus(id-s)P$ as $B$-modules;
cf.~Lemma~3(3) in~\cite{r}.\\

\bco\Label{galois} 
\slq\ is a faithfully flat \af\gal\ of \fr .
\eco 
\bpf
The fact that \slq\ is an \af\gal\ of \fr\ can be inferred from the  
proof of Proposition~\ref{coin}.\\ 

Another way to see it is as follows: For any Hopf algebra $P$, 
the canonical map 
$P\ot P\ni p\ot p'\to pp'\1\ot p'\2\in P\ot P$ is bijective. Consequently,
for any Hopf ideal $I$ of $P$, 
the canonical map $P\ot\sb{P\sp{co (P/I)}} P\ra P\ot (P\slash I)$ is surjective.
(Here we assume the natural right coaction 
$(id\te\pi)\ci\Delta :\, P\ra P\ot(P/I)\,$.)
Now, since in our case we additionally have that $P/I=\af$ is finite dimensional,
 we can conclude that \slq\ is an \af\gal\ of \fr\
by Proposition~\ref{coin} and \cite[Theorem~1.3]{s3} (see~\cite{kt}).\\

The faithful flatness of \slq\ over \fr\ follows from the commutativity of
the latter and Corollary~1.5 in \cite{s3} (see~\cite{kt}).
\epf

\bre\Label{ff}\em 
Note that just as the fact that \fr\ is the space of all coinvariants implies that
\slq\ is faithfully flat over it, the faithful flatness of \slq\ over \fr\ entails,
by virtue of \cite[Lemma~1.3(2)]{s1} (or the centrality of \fr\ in \slq\ and
 \cite[Remark~1.6(1)]{s2}), that \fr\ is the space of all coinvariants. 
Therefore it suffices
either to find all the coinvariants or prove the faithful flatness.
\ere

\bco\Label{strict} 
Sequence (\ref{frs}) is a strictly exact sequence
of Hopf algebras.
\eco

\bco\Label{cover} 
Sequence (\ref{frs}) allows one to view $SL_q(2)$ as a quantum group covering
of $SL(2,\IC )$ (see Section~18 in~\cite{pw}).
\eco

\bre\Label{local}\em 
We can think of $SL_q(2)$ as a quantum principal bundle over $SL(2,\IC)$.
This bundle, however, is {\em not} locally trivial 
in the sense of~\cite[p.460]{d}.
Indeed, otherwise it would have to be reducible to its classical subbundle
(see p.466 in \cite{d}), which is impossible because $SL_q(2)$ 
 has ``less" 
classical points (characters of \slq) than~$SL(2,\IC)$. 
(Cf. Section 4.2 in~\cite{bk}.)
\ere

\section{Quotients of 
\boldmath $A(SL\sb{e\sp{\mbox{\tiny $\frac{2\pi i}{3}$}}}(2))$ 
as cleft Hopf-Galois extensions}

Let us now consider the case of (quantum) Borel subgroups. To abbreviate
notation, in analogy with the previous section, we put 
$P\sb+=\slq\slash\lan c\ran\, ,\; B\sb+=\slc\slash\lan\bar{c}\ran$, and 
$H\sb+=P\sb+\slash\lan a\sp3-1,\, b\sp3\ran= A(F)\slash\lan\tilde{c}\ran$. 
(We abuse the notation by not distinguishing formally generators of
$P,\, P_+,\, P_-,\, P\sb\pm$, etc.)
As in the previous section, 
we have the Frobenius 
homomorphism
 (cf.~\cite[Section~7.5]{pw}) $Fr\sb+:B\sb+\ra P\sb+$ given by
the same formula as (\ref{fr}), and the associated exact sequence of Hopf
algebras:
\[
B\sb+\st{Fr\sb+}{\lra}P\sb+\st{\pi\sb+}{\lra}H\sb+\, .
\]
Before proceeding further, let us first establish a basis
of $P_+$ and a basis of~$H_+\,$.

\bpr \Label{+basis}
The set  $\{a^p b^r\}\sb{p,r\in \IZ,\, r\geq 0}$ is a 
basis of $P_+$.
\epr 
\bpf
This proof is based on the Diamond Lemma (Theorem~1.2 in~\cite{b}).
Let $\IC\lan\ha , \hb, \hd\ran$ be the free unital associative algebra generated
by $\ha , \hb, \hd$. We well-order the monomials of $\IC\lan\ha , \hb, \hd\ran$
first by their length, and then ``lexicographically" choosing the following 
order among letters: $\ha\preceq  \hd\preceq \hb $.
In particular, this is a semigroup partial ordering having descending chain
condition, as required by the Diamond Lemma. Furthermore, we chose the
reduction system ${\cal S}$ to be:
\[
{\cal S}=\left\{ (\ha\hd , 1)\, ,\, (\hd\ha , 1)\, ,\, 
(\hb\ha , q\sp{-1}\ha\hb)\, ,\, (\hb\hd , q\hd\hb)\right\}.
\]
It is straightforward to check that the aforementioned well-ordering is
compatible with ${\cal S}$, there are no inclusion ambiguities in  ${\cal S}$,
and all overlap ambiguities of ${\cal S}$ are resolvable. Therefore, by the
Diamond Lemma, the set of all ${\cal S}$-irreducible monomials 
is a basis of 
\mbox{$\IC\lan\ha , \hb, \hd\ran\slash J$}, 
$J:=\lan\ha\hd -1\, ,\,\hd\ha -1\, ,\,\hb\ha -q\sp{-1}\ha\hb\, ,
\,\hb\hd -q\hd\hb\ran$. The monomials 
$\ha\sp p\hb\sp r, \hd\sp k\hb\sp l$, \mbox{$ p,r,k,l\in \IN\sb0$}, $k>0$, are
irreducible under ${\cal S}$ and their image under the canonical surjection spans
\mbox{$\IC\lan\ha , \hb, \hd\ran\slash J$}. Consequently, they form a 
basis of $\IC\lan\ha , \hb, \hd\ran\slash J$. To conclude the proof it 
suffices to note that the algebras $\IC\lan\ha , \hb, \hd\ran\slash J$ and 
$P\sb+$ are isomorphic.
\epf

\bpr \Label{basis}
The set  $\{\tilde{a}^p \tilde{b}^r\}\sb{p,r\in\{0,1,2\}}$ is a basis of
$H_+$.
\epr 
\bpf
Analogous to the proof of Proposition~\ref{f-basis}.
\epf\ \\
 
The formula for the right coaction of $H_+$ on $P_+$
is not as complicated as (\ref{co}) and reads: 
 \beq \Label{coab}
\Delta_R (a^pb^r) = 
\sum_{\mu=0}^r  
\left({\scriptsize \ba{c}\!\!\! r\!\!\! \\ \!\!\! \mu\!\!\! \ea }\right)_q
q^{-\mu(2r -2 \mu)} a^{p+ \mu} b^{r-\mu}  
\otimes\tilde{a}^{2r + p -2 \mu}  \tilde{b}^{\mu}\, .
\eeq
With the above formula at hand, 
it is a matter of a straightforward calculation to prove that
$P_+$ is an $H_+$\gal\ of $Fr_+(B_+)$. In particular, we have
$P\sb+\sp{co H\sb+}=Fr\sb+ (B\sb+)$. Moreover, since $P_+$ 
is generated by a group-like and a skew-primitive element, 
it is
a pointed Hopf algebra. Consequently (see p.291 in~\cite{s1}), we obtain:

\bpr\Label{+}
$P\sb+$ is a cleft $H\sb+$\gal\ of $Fr\sb+ (B\sb+)$. 
\epr

Our next step is to
construct a family of cleaving maps for this extension.
To simplify the notation, for the rest of this paper we will identify
$B_+$ with its image under $Fr_+$.
First, we construct a family of unital convolution invertible
 $B_+$-linear maps $\Psi_\nu:P_+\ra B_+$, and then employ
Corollary~\ref{psico}. It is straightforward to verify that, for any
function $\nu:\{0,1,2\}\ra\IZ$ satisfying $\nu(0)=0$, the family 
$\{\Psi_\nu\}$ of $B_+$-homomorphisms given by the formula
\beq\Label{nu}
\Psi_\nu(a^pb^r):=\hd_{0r}a^{3\nu(p)},\; p,r\in\{0,1,2\},
\eeq
fulfills
the desired conditions. The convolution inverse of $\Psi_\nu$ is provided
by (see~\cite[p.47]{ad}) 
\[
\Psi_\nu^{-1}(a^pb^r):=\hd_{0r}a^{-3\nu(p)},\;\;\; p,r\in\{0,1,2\},\;\;\;
\Psi_\nu^{-1}(wt):=\Psi_\nu^{-1}(t)S(w),\;\;\; t\in P_+,\;\;\; w\in B_+\, .
\]
Consequently, 
\beq\Label{cmap}
\Phi_\nu:H_+\ra P_+,\;\Phi_\nu(\ta^p\tb^r)=\Psi_\nu^{-1}(a^{p+r})a^pb^r
=a^{-3([p+r]_1+\nu([p+r]_2))+p}b^r, 
\eeq
where $3[p+r]_1+[p+r]_2=p+r,\; 0\leq [p+r]_2<3$,
is a family of cleaving maps. In particular, we can choose $\nu(1)=0,\;\nu(2)=1$.
Then we have:
\bea \Label{Phi}
&&\Phi(1) =1~,~ \Phi(\tilde{a})=a~,~ \Phi(\tilde{a}^2)=a^{-1}~,~
\Phi(\tilde{b})=b~,~\Phi(\tilde{b}^2)=a^{-3}b^2~, \nonumber\\
&& \Phi(\tilde{a}\tilde{b})=a^{-2}b~,~
\Phi(\tilde{a}^2\tilde{b})=a^{-1}b~,~
\Phi(\tilde{a}\tilde{b}^2)=a^{-2}b^2~,~
\Phi(\tilde{a}^2\tilde{b}^2)=a^{-1}b^2~ .
\eea 
\bre\Label{gen}\em
Here we  rely on the fact that the monomials
$a^pb^r,\; p,r\in\{0,1,2\}$, form a $B_+$-basis of $P_+$.
As can be proven with the help of the linear basis 
$\{a^pb^r\}_{p,r\in\IZ,\,r\geq 0}\,$, the set 
$\{a^pb^r\}_{p,r\in\{0,\dots,n-1\}}\,$ is a $B_+$-basis of $P_+$
for any $n$-th primitive odd root of unity. Hence our construction
of a family of cleaving maps can be immediately generalised to an 
arbitrary primitive odd root of unity.
\ere
Let us now apply Lemma~\ref{cycle} to calculate explicitly
the cocycle $\sigma\sb\Phi:H\sb+\ot H\sb+\ra B\sb+\,$:  
\beq\Label{coceq}
\matrix{
\sigma\sb\Phi(\ta\ot\ta)=a\sp3, &
\sigma\sb\Phi(\ta^2\ot\ta^2)=a\sp{-3}, & 
\sigma\sb\Phi(\tb\ot\tb^2)=a\sp{-3}b^3, \cr 
\sigma\sb\Phi(\tb\ot\ta\tb^2)=q^2a\sp{-3}b^3, &
\sigma\sb\Phi(\tb\ot\ta^2\tb^2)=qb^3, &
\sigma\sb\Phi(\tb^2\ot\tb)=a\sp{-3}b^3, \cr
\sigma\sb\Phi(\tb^2\ot\ta\tb)=qa\sp{-6}b^3, & 
\sigma\sb\Phi(\tb^2\ot\ta^2\tb)=q^2a\sp{-3}b^3, &
\sigma\sb\Phi(\ta\tb\ot\tb^2)=a\sp{-6}b^3, \cr
\sigma\sb\Phi(\ta\tb\ot\ta\tb^2)=q^2a\sp{-3}b^3, &
\sigma\sb\Phi(\ta\tb\ot\ta^2\tb^2)=qa\sp{-3}b^3, &
\sigma\sb\Phi(\ta^2\tb\ot\tb^2)=a\sp{-3}b^3, \cr
\sigma\sb\Phi(\ta^2\tb\ot\ta\tb^2)=q^2a\sp{-3}b^3, &
\sigma\sb\Phi(\ta^2\tb\ot\ta^2\tb^2)=qa\sp{-3}b^3, &
\sigma\sb\Phi(\ta\tb^2\ot\tb)=a\sp{-3}b^3, \cr
\sigma\sb\Phi(\ta\tb^2\ot\ta\tb)=qa\sp{-3}b^3, &
\sigma\sb\Phi(\ta\tb^2\ot\ta^2\tb)=q^2a\sp{-3}b^3, &
\sigma\sb\Phi(\ta^2\tb^2\ot\tb)=b\sp3, \cr
\sigma\sb\Phi(\ta^2\tb^2\ot\ta\tb)=qa\sp{-3}b^3, &
\sigma\sb\Phi(\ta^2\tb^2\ot\ta^2\tb)=q^2a\sp{-3}b^3, &
\sigma\sb\Phi |\sb{\mbox{\scriptsize other basis elements}}=\he\ot\he\, .\cr
}
\eeq
As the cocycle action (see (\ref{act})) is necessarily trivial due to the
centrality of $B_+$ in $P_+\,$, we obtain the following: 

\bpr\Label{cp} 
$P_+$ is isomorphic as a comodule algebra to the twisted
product (see \cite[Example~4.10]{bcm}) of $B_+$ with $H_+$ defined
by the above described cocycle $\hs\sb\Phi\,$.
\epr

More explicitly, we can simply say that the algebra structure on $B_+\ot H_+$
that is equivalent to the algebra structure of $P_+$ is given by the formula
\beq\Label{mult}
(x\ot h)\cdot (y\ot l)=xy\hs\sb\Phi(h\1\te l\1)\ot h\2 l\2\, .
\eeq
Let us also mention that $\hs\sb\Phi$ is {\em not} a coboundary, i.e.,
it cannot be gauged by a unital convolution invertible map $\hg :H_+\ra B_+$
to the trivial cocycle $\he\ot\he$ (see Proposition~6.3.4 in~\cite{m-s}).
More formally, we have: 
\bpr\Label{coh}
The cocycle $\hs\sb\Phi$ represents a non-trivial cohomology class in the
(non-Abelian) 2-cohomology of $H_+$ with values in $B_+\,$.
\epr

\bpf
Suppose that the claim of the proposition is false. Then there would exist
\hg\ such that $\hs\sb{\gamma *\Phi}=\he\ot\he\,$, i.e., 
\beq\Label{stand}
[m\circ(\Phi\sp\gamma\te\Phi\sp\gamma)]*[(\Phi\sp\gamma)^{-1}\circ m]
=\he\ot\he\, .
\eeq
Here $\Phi\sp\gamma:=\gamma *\Phi$ and the middle convolution product is 
defined with respect to the natural coalgebra structure on $H_+\ot H_+$, namely
$\hD\sp\ot :=(id\ot\mbox{flip}\ot id)\circ(\hD\ot\hD)$. (Note that 
$\Phi\sp\gamma$ is also
a cleaving map.)
A standard argument (apply $*(\Phi\sp\gamma\circ m)$ from the right
to both sides of (\ref{stand})) allows us to conclude that $\Phi\sp\gamma$ is an
algebra homomorphism. Since $\Phi\sp\gamma$ is
injective (see Section~1), we can view 
$H_+$ as a subalgebra of $P_+\,$. In particular, there exists 
$0\neq p\in P_+$ such that $p^2=0$. (Put $p=\Phi\sp\gamma(\tb^2)$.)
Write $p$ as $\sum\sb{\mu\in{\IZ}} a\sp\mu p\sb\mu\,$, where the coefficients
$\{p\sb\mu\}\sb{\mu\in{\IZ}}$ are polynomials in $\tb$. 
Let $\mu_0(p):=\mbox{max}\{\mu\in{\IZ}\, |\; p\sb\mu\neq 0\}$.
It is well defined because $a\sp\mu b\sp n,\;\mu ,n\in{\IZ},\; n\geq 0$,
form a basis of $P_+\,$, and exists because $p\neq 0$. Now, due to the
commutation relation in $P_+$ and the fact that the polynomial ring 
${\IC}[b]$ has no zero divisors, we can conclude that $\mu_0(p^2)$ exists
(and equals $2\mu_0(p)$). This contradicts the equality $p^2=0$. 
\footnote{We are grateful to Ioannis Emmanouil for helping us to make the
nilpotent part of the proof simple.}
\epf\ \\

To put it simply, $H_+$ cannot be embedded in $P_+$ as a subalgebra.

\bre\em
Note that we could equally well try to use the {\it lower} (quantum) Borel 
subgroups $P_-$, $B_-$, $H_-$. The Hopf algebras $H_+$ and $H_-$ are naturally 
isomorphic as algebras and anti-isomorphic as coalgebras via the map 
that sends $\tilde{a}$ to $\tilde{a}$ and $\tilde{b}$ to $\tilde{c}$.
They are also isomorphic as coalgebras and anti-isomorphic 
as algebras via 
the map 
that sends $\tilde{a}$ to $\tilde{a}^2$ and $\tilde{b}$ to $\tilde{c}$. 
It might be worth noticing that $H_+$ and $H_-$ are {\em not}
 isomorphic as Hopf algebras. Indeed,
if they were so, there would exist an invertible algebra map $\phi 
:H_+ \ra H_-$ commuting with the antipodes. From direct computations, it 
turns out that any such map has to satisfy $\phi(\tilde{b})=
\kappa (\tilde{a} -q^2\tilde{a}^2)\tilde{c}^2 $, with $\kappa$ an arbitrary 
constant. This implies $\phi(\tilde{b})^2= \phi(\tilde{b}^2)=0$
 contradicting, due to $\tilde{b}^2 \neq 0$
(see Proposition~\ref{basis}), the injectivity of $\phi$ . 
\ere

To end this section, let us consider the Cartan case: We define the Hopf
algebras $P\sb\pm,\, B\sb\pm$ and $H\sb\pm$ by putting the off-diagonal
generators to 0, i.e., $P\sb\pm:=P/\lan b,c\ran,\,
 B\sb\pm:=B/\lan\bar{b},\bar{c}\ran,\, H\sb\pm:=H/\lan\tb,\tc\ran$.
Everything is now commutative, and we have $P_\pm\cong B_\pm\cong A(\IC\sp\times)$, 
$H_\pm\cong A(\IZ_3)$, where $\IC\sp\times :=\IC\setminus\{0\}$. 
It is immediate to see that, just as in the above discussed Borel case,
we have an exact sequence of Hopf algebras 
$B\sb\pm\st{Fr\sb\pm}{\ra}P\sb\pm\ra H\sb\pm\,$, and $P\sb\pm$ is a cleft
$H\sb\pm$\gal\ of $Fr\sb\pm(B\sb\pm)$. A cleaving map $\Phi$
and cocycle $\hs\sb\Phi$ are given by the formulas that look exactly as the
$a$-part of (\ref{Phi})
and (\ref{coceq}) respectively.
It might be worth to emphasize that, even though this extension is cleft, 
the principal  bundle 
$\IC\sp\times(\IC\sp\times , \IZ_3)$ is {\em not} trivial.
Otherwise $\IC\sp{\times}$ would have to be disconnected.
This is why we call $\Phi$ a cleaving map rather than a trivialisation.

\section{The bicrossproduct structure of \boldmath $P_+$}

In particular, the concept of cocleftness applies to the \hge s obtained
from short exact sequences of Hopf algebras. One can view cocleftness
as dual to cleftness the same way crossed coproducts are dual to
crossed products~\cite{m-s90}. 
The upper Borel extension $B_+\inc P_+\,$ is cleft and cocleft (see
Corollary~\ref{psico}), and
the maps $\Psi_\nu$ of (\ref{nu}) are both unital and counital. 
By Proposition~3.2.9 in~\cite{ad}, the cocleftness implies that
$P_+$ is isomorphic as a left $B_+$-module coalgebra to the crossed coproduct
of $B_+$ and $H_+$ given by the weak coaction
\[
\hl:H_+\ra H_+\ot B_+\, ,~~~ \hl(\pi_+(p)):=\pi_+(p\2)\ot \Psi^{-1}(p\1)\Psi(p\3)\, ,
\]
 and the co-cocycle 
\[
\hz:H_+\ra B_+\ot B_+\, ,~~~ 
\hz(\pi_+(p)):=\Delta(\Psi^{-1}(p\1))(\Psi(p\2)\ot\Psi(p\3))\, .
\]
Here $\Psi$ is the retraction obtained from (\ref{nu}) for the choice
of $\nu$ made above~(\ref{Phi}).
Explicitly, we have:  
\bea && \hl(\tilde{a}^p)=\tilde{a}^p\ot 1~,~
\hl(\tilde{b})=\tilde{b}\ot 1~,~ 
\hl(\tilde a \tilde{b})=\tilde a \tilde{b}\ot a^{-3}~,~ 
\hl(\tilde a ^2 \tilde{b})=\tilde a ^2\tilde{b}\ot a^{-3}~,\nonumber \\
&& \hl(\tilde{b}^2)=\tilde{b}^2\ot a^{-6}~,~ 
\hl(\tilde a \tilde{b}^2)=\tilde a \tilde{b}^2\ot  a^{-3}~,~ 
\hl(\tilde a ^2 \tilde{b}^2)=\tilde a ^2\tilde{b}^2\ot   a^{-3}~.
\eea
The co-cocycle is trivial, i.e., $\hz(\pi_+(p))=\he(p)\ot 1$.
We have thus arrived at:
\bpr 
$P_+$ is isomorphic as a left $B_+$-module coalgebra to the
 crossed coproduct  $B_+ \,^\lambda \# H_+$ defined by the above 
coaction $\hl$.
\epr
In particular, this means that the coproduct on $B_+\ot H_+$ 
that makes it isomorphic to $P_+$ as a coalgebra is
given by
\[
\Delta(w\ot h)=w\1\ot\hl^{[1]}(h\1)\ot w\2\hl^{[2]}(h\1)\ot h\2~,
\]
where $\hl(h):=\hl^{[1]}(h)\ot\hl^{[2]}(h)$ (summation suppressed).

\bpr\Label{corem}
The above defined coproduct is {\em not} equivalent to the tensor
coproduct  
$\hD^\otimes:=(id\ot\mbox{\em flip}\ot id)\ci(\hD\ot\hD): B_+\ot H_+\lra
(B_+\ot H_+)\ot (B_+\ot H_+)$.
\epr\bpf
Suppose the contrary. Then, by \cite[Proposition~3.2.12]{ad}, there would exist a
counital convolution invertible map $\xi:H_+\ra B_+$ such that
$\hl(h)= h\2 \ot\xi^{-1}(h\1)\xi(h\3)\,$.
With the help of Proposition~\ref{basis},
applying this formula to $\tilde b$ implies
$\xi^{-1}(\tilde a)\xi(\tilde a ^2)=1$, and requiring it for $\tilde b ^2 $
gives $\xi^{-1}(\tilde a ^2)\xi(\tilde a)=a ^{-6}$. 
Since $\ta$ is group-like, $\xi(\ta)$ and $\xi(\ta^2)$ are invertible,
and we obtain $1=\xi^{-1}(\tilde a)\xi(\tilde a ^2)=a^6$. This contradicts
Proposition~\ref{+basis}.
\epf\ \\

Note that as far as the algebra structure of $P_+$ is concerned, it is
given by the trivial action and a non-trivial cocycle. For the coalgebra
structure it is the other way round, i.e., it is given by the trivial
co-cocycle and a non-trivial coaction (cf.~\cite{m-s97}). 
Due to the triviality of co-cocycle $\hz$, $\,\Psi$ is a coalgebra
homomorphism. Also, one can check that 
the cocycle and coaction put together make $B_+\ot H_+$  
a {\em cocycle bicrossproduct Hopf algebra} \cite{m-s}. 

\bco\Label{bico}
The Hopf algebra $P_+$ is isomorphic to the cocycle bicrossproduct Hopf
algebra  $B_+\,^\lambda \#_\sigma H_+\,$. The  
isomorphism and its inverse are given by 
\[ 
p\mapsto \Psi(p\1)\ot \pi_+(p\2)~,~~~w\ot h\mapsto w\Phi(h)~.
\] 
Here $\Phi$ is related to $\Psi$ as in Corollary~\ref{psico}, and
given explicitly by formulas~(\ref{Phi}).
\eco

\section{Integrals on and in \boldmath $A(F)$}

Recall that a left (respectively right) integral
{\em on} a Hopf algebra $H$ over a field $k$ 
is a linear functional $h:H \ra k$ 
satisfying:
\beq\Label{int}
(id \ot h) \circ \Delta = 1\sb H \cdot h ~~~
(\mbox{respectively~} (h\ot id)\circ \Delta = 1\sb H \cdot h ).
\eeq 
(For a comprehensive
review of the theory of integrals see \cite[Section~1.7]{m-s}, 
\cite[Chapter~V]{s}.)
In the case of the Hopf algebra \af , we have the following result:
\bpr
The space of left integrals on \af\ coincides with the space of right 
integrals on \af , i.e., \af\ is a unimodular Hopf algebra. 
In terms of the basis $\{\ta^p\tb^r\tc^s\}_{p,r,s\in\{ 0,1,2\}}$ of \af ,
for any integral $h$, we have by 
$h(\tilde{a}^p \tilde{b}^r\tilde{c}^s)=z\hd^p_0\hd^r_2\hd^s_2\, ,\; z\in\IC$. 
\epr 
\bpf 
By applying the projection $\pi_{\pm}:\af\ra H_{\pm}$ to (\ref{int}), it is 
easy to see that any left (and similarly any right) integral has to vanish 
on about half of the elements of the basis. With this information at hand,
and using the fact that on a finite dimensional Hopf algebra the space of 
left and the space of right integrals are one dimensional \cite{ls},
it is straightforward to verify by a direct calculation the claim of the
proposition. 
\epf\ \\

A two-sided integral on a Hopf algebra $H$ is called a {\it Haar measure} 
iff it is {\it normalised}, i.e., iff  $h (1)= 1$.
As integrals on \af\ are {\em not} normalisable, we have:
\bco
There is no Haar measure on the Hopf algebra \af\
(cf.~Theorem~2.16 in~\cite{kp} and (3.2) in~\cite{mmnnu}).
\eco

\bre\em 
Since the Hopf algebra \af\ is finite dimensional, $F$ can be considered
as a finite quantum group. However, it is {\em not} a compact matrix 
quantum group in the sense of Definition~1.1 in~\cite{w-s}. Indeed, by
Theorem 4.2 in~\cite{w-s},
compact matrix quantum groups always admit a (unique) Haar measure.
Furthermore, as \af\ satisfies all the axioms of Definition~1.1 in~\cite{w-s}
except for the $C^*$-axiom, there does not exist 
a $*$-structure and a norm on \af\  that would make \af\  a 
Hopf-$C^*$-algebra. 
In particular, for the $*$-structure given by setting 
$\tilde{a}^*=\tilde{a}$, 
$\tilde{b}^*=\tilde{b}$, 
$\tilde{c}^*=\tilde{c}$, $\tilde{d}^*=\tilde{d}$, this fact is evident: 
Suppose that there exists a norm 
satisfying the $C^*$-conditions. Then $0=\parallel \tilde{c}^4 
\parallel =
\parallel (\tilde{c}^2)^*\tilde{c}^2 \parallel = \parallel \tilde{c}^2 
\parallel^2$, which implies $\tilde{c}^2=0$ and thus contradicts 
Proposition~\ref{f-basis}. 
\ere 
 
We recall also that an element $\Lambda \in H$ is called a left 
(respectively right) integral 
{\it in} $H$, iff it verifies $\alpha\Lambda =\varepsilon(\alpha)\Lambda$,
(respectively 
$\Lambda\alpha =\varepsilon(\alpha)\Lambda$) for any $\alpha \in H$.
If $H$ is finite dimensional,  an integral in $H$ corresponds to 
an integral on the dual Hopf algebra $H^*$. Clearly, an integral in \af\
should annihilate any non-constant polynomial in 
$\tilde{b}$ and $\tilde{c}$, whereas it should leave unchanged any 
polynomial in $\tilde{a}$. It is easy to see that the element 
$\Lambda_L =(1+\tilde{a}+\tilde{a}^2)
\tilde{b}^2\tilde{c}^2$ is a left integral and the element  
$\Lambda_R =
\tilde{b}^2\tilde{c}^2 (1+\tilde{a}+\tilde{a}^2)$ is a right integral. 
Hence in this case left and right integrals are not proportional. We
can therefore
conclude that $H^*$, which by Section~3 in~\cite{dns} can be identified
with $U\sb q(sl\sb2)/\lan K^3-1,E^3,F^3\ran$ of~\cite{c-r}, 
is {\em not} unimodular. 
Again, since
\af\ is finite dimensional, any left integral in \af\ is proportional to
$\Lambda_L\,$, and any right integral in \af\ is proportional to $\Lambda_R\,$.
In addition, by Theorem 5.1.8 in~\cite{s},
the property $\varepsilon(\Lambda_L)=0$ assures us that \af\ is {\it not} 
semisimple as an algebra. \\

\section{A coaction of \boldmath $A(F)$ on \boldmath $M(3,\IC)$}

Let us now consider $F$ as a quantum-group symmetry of $M(3,\IC)$
--- a direct summand of A.~Connes' algebra for the Standard Model.
Recall first that for any $n\in\IN$ the algebra of matrices $M(n,\IC)$
can be identified with the algebra
\mbox{$\IC\lan\mbox{x},\mbox{y}\ran\slash 
\lan\mbox{xy}-\mu\mbox{yx},\,\mbox{x}\sp n-1,\,\mbox{y}\sp n-1\ran$},
$\mu=e\sp{\frac{2\pi i}{n}}$. (Map x to 
${\scriptsize\left(\ba{cc}0 & I\sb{n-1}\\ 1 & 0\ea\right)}$
and y to $diag(1,\mu,...,\mu\sp{n-1})$; see Section~IV.D.15 of~\cite{w-h}.)
Denoting by ${A}(\IC_q^2)$ the polynomial algebra (in $x$ and $y$)
of the quantum plane, by 
${A}(\IC^2)=\IC[\bar{x},\bar{y}]$ the algebra of polynomials on $\IC\sp2$,
and maintaining our assumption that $q=e\sp{\frac{2\pi i}{3}}$, we obtain
the following sequence of algebras and algebra homomorphisms: 
\beq\Label{fra}
 {A}(\IC^2) \mathop{\lra}\limits^{fr} {A}(\IC_q^2)
\mathop{\lra}\limits^{\pi\sb M} M(3,\IC)\cong{A}(\IC_q^2)\slash
\lan x\sp3-1,\, y\sp3-1\ran  \ .
\eeq
Here $fr$ is an injection given by $fr(\bar{x})=x\sp3,\; fr(\bar{y})=y\sp3$,
and $\pi\sb M$ is the map induced by the canonical surjection
\mbox{${A}(\IC_q^2)\ra {A}(\IC_q^2)\slash
\lan x\sp3-1,\, y\sp3-1\ran$}.\\

Let us note that although ${A}(\IC_q^2)$ and ${A}(\IC^2)\otimes M(3,\IC)$
are isomorphic as ${A}(\IC^2)$-modules, their algebraic structures 
(cf.~(\ref{mult}))
are slightly different:
\beq\Label{plmul}
 (\bar x^p \bar y^r\ot\tx^k\ty^{\ell})
(\bar x^s \bar y^t\ot\tx^m\ty^n ) 
 =\bar x^{p+s+[k+m]_1} \bar y^{r+t+[\ell +n]_1}
\ot\tx^{[k +m]_2}\ty^{[\ell +n]_2}\, ,
\eeq
where $3[n]_1+[n]_2=n,\; [n]_1,[n]_2\in\IN,\; 0\leq [n]_2<3,\;
\tx=\pi_M(x),\; \ty=\pi_M(y)$.
Incidentally, the associativity of this product amounts to the identity
$[k\! +\! m]_1 \! +\! [[k\!  +\! m]_2\!  +\! u]_1 = 
[m\! +\! u]_1\! +\! [[k\! +\! [m\! +\! u]_2\! +\! u]_1$.\\

Next, observe that combining sequences (\ref{fra}) and  (\ref{frs})
 together with the natural right coactions 
($\mbox{e}\sb i\to\sum\sb{j\in\{1,2\}}\mbox{e}\sb j\ot M\sb{ji}\, ,\; i\in\{1,2\}$)
on  ${A}(\IC^2)$,  ${A}(\IC_q^2)$, and $M(3,\IC)$ respectively, one can
obtain the following commutative diagram of algebras and algebra homomorphisms:
\beq
\begin{array}{ll}
 {A}(\IC^2) \!\!\! &\!\!\!\! \mathop{-\hspace{-6pt}
\longrightarrow}\limits^{\rho} {A}(\IC^2) \otimes {A}(SL(2,\IC))  \\
\ ^{fr} \Big\downarrow 
 &  ~~~~~~~~~~~~~~~\Big\downarrow \phantom{\cdot} ^{fr\otimes F\!r} \\ 
\ {A}(\IC_q^2) \!\!\! &\!\!\!\! \mathop{-\hspace{-6pt}
\longrightarrow}\limits^{\rho_q} 
{A}(\IC^2_q) \otimes {A}(SL_q(2,\IC)) \\
\ ^{\pi\sb M} \Big\downarrow ~
 &  ~~~~~~~~~~~~~~~\Big\downarrow \phantom{\cdot} ^{\pi\sb M\otimes 
\pi\sb F} \\ \ M(3,\IC) \! &\!\!\!\! \mathop{-\hspace{-6pt}
\longrightarrow}\limits^{\rho_F} M(3,\IC) \otimes {A}(F)  \ .
\end{array}
\eeq\ \\

Another way to look at $M(3,\IC)\st{\rho_F}{\ra}M(3,\IC)\ot\af$ is to treat
$M(3,\IC)$ as a 9-dimensional comodule rather than a comodule algebra. Let
us choose the following linear basis of $M(3,\IC)$:
\[
e_1=1,\; e_2=\tx,\; e_3=\ty,\; e_4=\tx^2,\; e_5=\tx\ty,\; e_6=\ty^2,\; 
e_7=\tx^2\ty,\; e_8=\tx\ty^2,\; e_9=\tx^2\ty^2.
\]
The formula $\dr e_i=e_j\te N_{ji}$
allows us to determine the corepresentation matrix $N$:
\beq\Label{n}{\footnotesize
\pmatrix{
1 & 0 & 0 & 0 & 0 & 0 & \ta^2(\tb+q^2\tc^2) & \ta(\tb^2+q^2\tc-q\tb\tc^2) & 0 \cr
0 & \ta & \tb & 0 & 0 & 0 & 0 & 0 & \ta^2(\tb^2-q\tc) \cr
0 & \tc & \td & 0 & 0 & 0 & 0 & 0 & \ta(q^2\tb^2\tc+q\tc^2-\tb) \cr
0 & 0 & 0 & \ta^2 & \ta\tb & \tb^2 & 0 & 0 & 0 \cr
0 & 0 & 0 & -q^2\ta\tc & (1-\tb\tc) & -q^2\tb\td & 0 & 0 & 0 \cr
0 & 0 & 0 & \tc^2 & \tc\td & \td^2 & 0 & 0 & 0 \cr
0 & 0 & 0 & 0 & 0 & 0 & \ta & -\tb & 0 \cr
0 & 0 & 0 & 0 & 0 & 0 & -\tc & \td & 0 \cr
0 & 0 & 0 & 0 & 0 & 0 & 0 & 0 & 1 \cr
}}
\eeq
It is clear that $N$ is reducible. The upper right corner terms of $N$ appear
to be an effect of the finiteness of $F$.  By restricting the comodule $M(3,\IC)$
respectively to the linear span of $1,\tx^2\ty,\tx\ty^2$  and 
the linear span of
$\tx,\ty,\tx^2\ty^2$, we obtain two ``exotic"
corepresentations of~\af\ (see~\cite[Section~4]{dns} for the dual picture):
\beq\Label{exotic}
N_1=
\pmatrix{
1 & \ta^2(\tb+q^2\tc^2) & \ta(\tb^2+q^2\tc-q\tb\tc^2) \cr
0 & \ta & -\tb \cr
0 & -\tc & \td \cr
},\;\;\;
N_2=
\pmatrix{
\ta & \tb & \ta^2(\tb^2-q\tc) \cr
\tc & \td & \ta(q^2\tb^2\tc+q\tc^2-\tb) \cr
0 & 0 & 1 \cr
}.
\eeq\ \\

To end with, let us remark that, very much like the Frobenius map $Fr$,
the ``Frobenius-like" map $fr$ of sequence (\ref{fra}) allows us to identify
$A(\IC\sp2)$ with the subalgebra of 
$(id\ot\pi\sb F)\circ\rho\sb q$-coinvariants of $A(\IC\sp2\sb q)$:
\beq\Label{eq}
fr(A(\IC\sp2))=A(\IC\sp2\sb q)\sp{co A(F)}.
\eeq
Indeed, since we can embed $A(\IC\sp2\sb q)$ in \slq\ as a subcomodule
algebra (e.g., \mbox{$x\to a$}, \mbox{$y\to b$}), equality (\ref{eq})
 follows directly from Proposition~\ref{coin} and the lemma below:

\ble
Let $P\sb1$ and $P\sb2$ be right $H$-comodules, and 
\mbox{\em j $:P\sb1\ra P\sb2$} an injective comodule homomorphism.
Then $P\sb1\sp{co H} = \mbox{\em j}\sp{-1}(P\sb2\sp{co H})$.
\ele\bpf
Denote by \mbox{$\rho\sb1:P\sb1\ra P\sb1\ot H$} and
\mbox{$\rho\sb2:P\sb2\ra P\sb2\ot H$} the right $H$-coactions on 
$P\sb1$ and $P\sb2$ respectively. Assume now that $p\in P\sb1\sp{co H}$.
Then $\rho\sb2(\mbox{j}(p))=(\mbox{j}\ot id)(\rho\sb1(p))= \mbox{j}(p)\ot 1$,
i.e., $p\in\mbox{j}\sp{-1}(P\sb2\sp{co H})$. Conversely, assume that
$p\in\mbox{j}\sp{-1}(P\sb2\sp{co H})$. Then  
$(\mbox{j}\ot id)(p\ot 1)=\rho\sb2(\mbox{j}(p))=(\mbox{j}\ot id)(\rho\sb1(p))$.
Consequently, by the injectivity of $(\mbox{j}\ot id)$, we have 
$\rho\sb1(p)=p\ot 1$, i.e., $p\in P\sb1\sp{co H}$.
\epf\\

{\bf Acknowledgments:} 
P.M.H. was partially supported by a visiting fellowship at 
the International Centre for Theoretical Physics in Trieste, 
KBN grant \mbox{2 P301 020 07} and NATO postdoctoral fellowship.
It is a pleasure to thank Tomasz Brzezi\`nski, 
Peter Schauenburg and Hans-J\"urgen Schneider
for crucial discussions. 
We are also grateful to Ioannis Emmanouil and Shahn Majid for their advice.

\end{document}